\input amstex
\documentstyle{amsppt}
\magnification=1200

\catcode`\@=12


\def\ra{\rightarrow}
\def\kc{{\Cal K}}
\def\lc{{{\Cal L}}}

\def\oc{{\Cal O}}
\def\fg{{\frak g}}

\def\fv{{\frak v}}

\def\ei{\epsilon_{1}}
\def\e2{\epsilon_{2}}
\def\ez{\epsilon_{2}}

\def\pa{\pi_{\la_{*}}}

\def\bp{b^\prime}

\def\l2{L^2}
\def\rb{{\bold R}}

\def\trace{\text{trace}}
\def\g.c.d.{\text{g.c.d.}}
\def\ker{\text{ker}}

\def\so{\text{so}}
\def\det{\text{det}}
\def\trace{\text{trace}}
\def\span{\text{span}}
\def\spec{\text{spec}}
\def\Vol{\text{vol}}
\def\Ric{\text{Ric}}

\def\cl{{{\Cal L}}}
\def\jp{j^\prime}
\def\jo{j_0}

\def\ju{j_u}
\def\ogj{\overline{G(j)}}
\def\ogjp{\overline{G(j^\prime)}}
\def\mj{M(j)}
\def\mjp{M(j^\prime)}
\def\mjl{M_\lambda(j)}
\def\mlj{M_\lambda(j)}
\def\mjlp{M_{\lambda}(j^{\prime})}

\def\gj{G(j)}
\def\gjp{G(j^\prime)}
\def\oglj{\overline{G_\lambda(j)}}
\def\gljp{G_\lambda(j^{\prime})}
\def\glj{G_\lambda(j)}
\def\ogljp{\overline{G_\lambda(j^{\prime})}}
\def\fz{{\frak{z}}}

\def\la{\lambda}

\def\cx{{\Cal X}}

\def\tp{T^\prime}
\def\fgl{{\frak g\frak l}}

\def\fgs{{\frak g}}
\def\gg{{\Gamma\backslash G}}
\def\bcb{<\,,\,>}

\def\ll{{\lambda}}
\def\la{\lambda}
\def\fzl{\fz_{\ll}}
\def\Hl{{\Cal H}_\ll}
\def\hl{{\Cal H}_{[\ll]}}
\def\Hlp{{\Cal H}_\ll^\prime}
\def\hlp{{\Cal H}_{[\ll]}^\prime}
\def\njg{N(j,\Gamma)}
\def\njp{N(j^\prime,\Gamma^\prime)}
\def\nlj{N_\ll(j,\Gamma)}
\def\nljp{N_\ll(j^\prime,\Gamma^\prime)}

\newcount\theoremnumber
\def\grow{\advance \theoremnumber by 1}
\def\theorem#1\par{\grow
\noindent {\bf Theorem}\quad {\the\theoremnumber} .\qquad {\sl #1}\par}

\topmatter

\title Continuous families of isospectral Riemannian metrics which are not
locally
isometric \endtitle

\rightheadtext{Isospectral Riemannian metrics}

\author Carolyn S. Gordon and Edward N. Wilson
\endauthor

\affil Dartmouth College and Washington University
\endaffil

\address 
\flushpar Carolyn S. Gordon:  Dartmouth College,  Hanover, New Hampshire, \
03755;
\newline carolyn.s.gordon@dartmouth.edu
\flushpar Edward N. Wilson: Washington University, St. Louis, Missouri\ 63130;
\newline enwilson@math.wustl.edu
\endaddress

\keywords Spectral Geometry, isospectral deformations, nilpotent Lie groups
\endkeywords
\thanks The first author is partially supported by a grant from
the National Science Foundation.  
\newline This research was initiated at the Mathematical
Sciences Research Institute.  Research at MSRI is supported by NSF grant \#DMS
9022140.
\endthanks
\endtopmatter

\document
\subheading{Introduction}
\medskip  

Two compact Riemannian manifolds are said to be {\it isospectral} if
the associated Laplace-Beltrami operators, acting on smooth functions, have the
same
eigenvalue spectrum.  If the manifolds have boundary, we specify {\it
Dirichlet} or
{\it Neumann} isospectrality depending on the boundary conditions imposed on
the
eigenfunctions.   

Numerous examples of isospectral compact manifolds have been constructed; see,
for
example, \cite{BT}, \cite{Bu}, \cite{DG}, \cite{GW 1,2}, \cite{GWW},
\cite{Gt1,2}, \cite{I}, \cite{M}, \cite{S} and \cite{V} or the survey articles
\cite{Be}, \cite{B}, \cite{D}, and \cite{G1}. Until recently however, all known
examples of isospectral manifolds were locally isometric, though not globally
isometric. In particular, the closed isospectral manifolds had a common cover.
Then Z.
Szabo \cite{Sz} gave a construction of pairs of isospectral compact manifolds
with
boundary which are not locally isometric, and the first author \cite{G2,3}
constructed pairs of isospectral closed Riemannian manifolds which are not 
locally
isometric.   Szabo pointed out that the curvature operators of these
isospectral
manifolds have different eigenvalues, thus identifying a specific local
invariant
which is not spectrally determined.

The first goal of this paper is to exhibit {\it continuous} families of
isospectral
Riemannian manifolds which are not locally isometric, i.e., we continuously
deform the
Riemannian metrics in such a way that the local geometry changes but the
Laplace
spectrum remains invariant. In
fact we prove:

\proclaim{Theorem 0.1}  Let $B$ be a ball in $\rb^{m}, m \ge 5$, and let
$T^{r}$ be
a torus of dimension $r \ge 2$. There exist continuous $d$-parameter families
of
Riemannian metrics on the compact manifold $B \times T^{r}$ which are both
Dirichlet
and Neumann isospectral but not locally isometric.  Here $d$ is
of order at least $O(m^2)$.
\endproclaim

A precise lower bound on $d$ is given in Theorem 2.2.

We also consider closed manifolds. Here we have not been able to construct
examples 
of continuous isospectral deformations in which the metrics are not locally
isometric.
However,  we do construct new examples of pairs of isospectral closed manifolds
which
are not  locally isometric. 

Next we examine the local geometry of the isospectral manifolds.  Since all the
manifolds considered in this paper are locally homogeneous, the curvature does
not vary from point to
point. In particular the eigenvalues of the Ricci tensor are constant functions
on each
manifold.  We exhibit specific examples of isospectral deformations of
manifolds with
boundary for which the eigenvalues of the Ricci tensor deform non-trivially. 
Similarly, we exhibit pairs of isospectral closed manifolds whose Ricci tensors
have
different eigenvalues.  These examples illustrate for the first time that the
Ricci
curvature is not spectrally determined. 

The paper is organized as follows:

In \S 1, we give a method for constructing isospectral metrics on $B \times
T^{r}$
which are not locally isometric. The construction reduces to a problem in
linear
algebra: 
\smallskip
\item{{\bf (P)}} Find pairs of $r$-dimensional subspaces of 
$so(m)$
and an isomorphism between these subspaces such that
corresponding elements have the same spectrum but such that the two
subspaces are not conjugate by any orthogonal transformation.
\smallskip

   As will be explained in \S1, each subspace of $so(m)$ gives rise to a
two-step
nilpotent Lie algebra with an inner product and thus to a simply-connected
nilpotent
Lie group with a left-invariant Riemannian metric.  The non-conjugacy condition
in (P)
guarantees that the resulting pair of nilpotent Lie groups with metrics are not
locally
isometric.  The manifolds in Theorem 0.1 are domains with boundary in these
nilpotent
Lie groups (more precisely, in nilpotent Lie groups covered by these
simply-connected
ones). We show that the spectral condition in (P) guarantees the isospectrality
of
these compact manifolds with boundary.  We end
\S 1 with a 7-dimensional example.

In \S 2, we give an explicit construction of continuous families of 
2-dimensional
subspaces of $\so(6)$ satisfying pairwise the condition (P) described above. 
Moreover, we show that for $m=5$ and for $m\geq 7$, generic two-dimensional
subspaces
of $so(m)$ belong to $d$-parameter families of subspaces which satisfy pairwise
the
condition (P), where $d$ is of order $O(m^2)$.  This completes the proof of
Theorem
0.1. 

In \S 3 we consider {\it nilmanifolds}, i.e., closed manifolds arising as
quotients
$\Gamma
\backslash G$ of nilpotent  Lie groups by discrete subgroups, endowed with
Riemannian metrics induced from left-invariant metrics on $G$.  We generalize
the
construction given in 
\cite{G2,3} of isospectral nilmanifolds. We construct seven and 
eight-dimensional
examples of isospectral nilmanifolds by  taking quotients of
suitable pairs of the simply-connected  nilpotent Lie groups occurring in the
examples
in \S 1 and \S 2.

\S 4 examines the curvature of the various examples, in particular showing 
that many
of the isospectral manifolds have different Ricci curvature. 

An appendix supplies a proof of a result needed in \S 3.

We wish to acknowledge Zoltan Szabo's beautiful work \cite{Sz} which inspired
Theorem
0.1.

\head \S1 Lie Algebra Criteria for Local Isometry and Isospectrality \endhead
 
A left-invariant Riemannian metric $g$ on a connected Lie group $G$ corresponds
to an inner product $<\cdot\,,\,\cdot>$ on the Lie algebra ${\frak g}$ of $G$.
We
will call the pair $({\frak g}, <\cdot\,,\,\cdot>)$ a {\it metric Lie
algebra}. Recall that
$G$ is said to be {\it two-step nilpotent} if $[{{\frak g}}, {\frak g}]$ is a
non-zero
subspace of the center of ${\frak g}$. Letting ${\frak z} = [{\frak g},
{\frak g}]$ and
${\frak v} = {\frak z}^{\perp}$
relative to $<\cdot\,,\,\cdot>$, we can then define an injective linear map
$j : \fz \ra \so({\frak v}, <\cdot\,,\,\cdot>)$ by
$$
\langle j(z) x, y \rangle = \langle [x, y], z \rangle\,\text{ for }\,  
x, \ y \in {\frak v},  z \in {\frak z}. \tag1.1
$$

Conversely, given any two finite dimensional real inner product spaces
${\frak v}$ and ${\frak z}$ along with a 
linear map $j : {\frak z} \ra \so({\frak v})$, we can define a metric Lie
alegbra
${\frak g}$ as the inner product space direct sum of ${\frak v}$
and ${\frak z}$ with the alternating bilinear bracket
map 
$[\cdot, \cdot] : {\frak g} \times {\frak g} \ra {\frak z}$ defined by
insisting that 
${\frak z}$ be central in ${\frak g}$ and using (1.1) to define $[x, y]$  for
$x, y \in {\frak v}$. Then ${\frak g}$ is two-step nilpotent if $j$ is non-zero
and
${\frak z} = [{\frak g},{\frak g}]$ if $j$ is injective.  We will always assume
$j$
is injective.

In the sequel, we will fix finite dimensional inner product spaces
${\frak v}$
and ${\frak z}$, use $<\,,\,>$ as a generic symbol for the fixed inner products 
on ${\frak v}, {\frak z}$ and ${\frak g} = {\frak v} \oplus {\frak z}$, and
we will contrast  properties of objects arising from pairs $j,\,j^{\prime}$ of
linear
maps from
${\frak z}$ to $\so({\frak v})$. 

\smallskip
 
\definition{Notation 1.2} 
(i) The metric Lie algebra defined as above from the data $({\frak v}, {\frak
z}, j)$ will be denoted ${\frak g}(j)$ and the corresponding simply-connected
Lie group
will be
denoted $G(j)$.  The Lie group $G(j)$ is endowed with the left-invariant
Riemannian
metric
$g$ determined by the inner product $<\cdot\,,\,\cdot>$ on ${\frak g}(j)$.

(ii) Explicitly, $G(j)$ may be identified diffeomorphically (though not
isometrically) with the Euclidean space ${\frak v} \times {\frak z}$ consisting
of all
pairs $(x, z)$ with $x\in {\frak v}, z \in {\frak z}$. The group product on
$G(j)$ is
given by
$$
(x,z)(x^{\prime}, z^{\prime}) = (x + x^{\prime}, z + z^{\prime} +
\frac{1}{2}[x, x^{\prime}])
$$
The Lie algebra element in ${\frak g}(j)$ determined by $x\in {\frak v},  z \in 
{\frak z}$ will be denoted by $x + z$ with the diffeomorphism $\exp: {\frak
g}(j) \ra
G(j)$ thereby  expressed by $\exp(x +z) = (x,z)$.  The exponential map
restricts to a
linear isomorphism between $\fz \subset \fg(j)$ and the derived group $[G,G]$
of $G$.

(iii) Suppose ${\Cal L}$ is a lattice of full rank in  ${\frak z}$, i.e. 
$\overline{{\frak z}} = {\frak z}/{\Cal L}$ is a torus. Denote by
$\overline{G(j)}$
the quotient of the Lie group $G(j)$ by the discrete central 
subgroup $\exp ({\Cal L})$. Then $\overline{G(j)}$ is again a connected Lie
group with
Lie algebra ${\frak g}(j)$. Diffeomorphically, $\overline{G(j)}$ may be
identified with
${\frak v}\times\overline{{\frak z}}$ and the exponential map 
$\overline{\exp} : {\frak g} (j) \ra \overline{G(j)}$ is expressed by 
$$
\overline{\exp}(x + z) = (x, \overline{z}) \text{ for } x \in v, z \in {\frak
z},
\text{ and } \overline{z} = z + {\Cal L} \in \overline{{\frak z}}.
$$
We assign to $\overline G(j)$ the unique left-invariant Riemannian metric
determined
by 
$<~\cdot\,,~\,\cdot~>$.  Thus the canonical projection from $G(j)$ to
$\overline G(j)$
given by
$(x,z) \to (x,
\overline{z})$ is  a Riemannian covering map as well as a Lie group
homomorphism.

(iv) For $B = \{x \in {\frak v} : \|x\| \leq 1\}$ the unit ball around 0
in ${\frak v}$ and for ${\Cal L}$ as in (iii), denote by $M(j)$ the subset 
$B \times \overline{{\frak z}} = \overline{\exp}(B+{\frak z})$ of
$\overline{G(j)}$ equipped with the inherited Riemannian structure.
$M(j)$ is thus a compact Riemannian submanifold of $\overline{G(j)}$ of full
dimension with boundary diffeomorphic to 
$S\times \overline{{\frak z}}$ for $S$ the unit sphere
around
$0$ in
${\frak v}$.  (Here we are using the identifications described in (iii). 
$M(j)$ of
course depends on the choice of $\lc$, but we view this choice as fixed.)
\enddefinition

\definition{Definition 1.3}  Let ${\frak v}$ and ${\frak z}$ be as above.  

(i) A pair $j,
j^{\prime}$ of linear maps from ${\frak z}$ to $\so({\frak v})$ will be called
{\it
equivalent}, denoted
$j \simeq j^{\prime}$, if there exist orthogonal linear operators $A$ on
${\frak v}$ 
and $C$ on ${\frak z}$ such that
$$
Aj(z)A^{-1} = j^{\prime}(C(z))
$$
for all $z \in {\frak z}$.

(ii) We will say $j$ is {\it isospectral} to $j^{\prime}$, denoted $j \sim
j^{\prime}$, if for each $z \in {\frak z}$, the eigenvalue spectra (with
multiplicities) of $j(z)$ and $j^{\prime}(z)$ coincide, i.e., there exists an
orthogonal linear operator
$A_z$ for which 
$$ {A_z j(z) A_z^{-1} = j^{\prime}(z)}.
$$
\enddefinition

\proclaim{Proposition 1.4}  Let ${\frak v}$ and ${\frak z}$ be finite
dimensional
real inner product spaces, $j$ and $j^{\prime}$ linear injections from ${\frak
z}$ to
$\so({\frak v})$, and ${\Cal L}$ a lattice of full rank in ${\frak z}$. Let 
${\frak g}(j)$, $G(j)$, $\overline{G(j)}$, and $M(j)$ be the objects defined in
1.2
from the data $({\frak v}, {\frak z}, j, {\Cal L})$ and let ${\frak
g}(j^{\prime})$,
$G(j^{\prime})$,
$\overline{G(j^{\prime})}$, and $M(j^{\prime})$ be the corresponding objects
defined by
the  data $({\frak v}, {\frak z}, \ j^{\prime}, {\Cal L})$. Then the following
are
equivalent:

\item{(a)} $\overline{G(j)}$ is locally isometric to
$\overline{G}(j^\prime)$;
\item{(b)} $M(j)$ is locally isometric to $M(j^{\prime})$;
\item{(c)} $G(j)$ is isometric to $G(j^{\prime})$;
\item{(d)} $j \simeq j^{\prime}$ in the sense of Definition
1.3.

\endproclaim

\demo{Proof} The local geometries of $G(j), \overline{G}(j)$, and $M(j)$
are identical. Thus each of $(a)$ and $(b)$ is equivalent to saying that $G(j)$
is
locally isometric to
$G(j^{\prime})$ which, by simple-connectivity, is equivalent to $(c)$. The
second
author showed in \cite{W} that if $(G,g)$ and $(G^{\prime}, g^{\prime})$ are
two
simply-connected nilpotent Lie groups with left-invariant metrics $g,
g^{\prime}$ and 
associated metric Lie algebras $({\frak g} <\cdot\,,\,\cdot>)$, $({\frak
g}^{\prime},
<\cdot\,,\,\cdot>^{\prime})$, then $(G,g)$ is isometric to $(G^{\prime},
g^{\prime})$,
if and only if there exist a map $\tau : {\frak g} \rightarrow {\frak
g}^{\prime}$
which is both a Lie algebra isomorphism and an inner product space isometry. In
our
case, equivalence of $(c)$ and $(d)$ follows by routine use of (1.1) serving to
reduce
these conditions on
$\tau$ to
$j \simeq j^{\prime}$.
\enddemo

\proclaim{Theorem 1.5}
Let ${\frak v}$ and ${\frak z}$ be inner product spaces, $j, \jp :{\frak z}\ra
\so(
{\frak v})$ linear injections, $\cl$ a lattice of full rank in ${\frak z}$, and
~
$\mj$~ and ~$\mjp$~ the manifolds defined in 1.2 from the data $({\frak
v},{\frak z},j,\cl)$ and
$({\frak v},{\frak z},j^\prime,\cl)$, respectively.  Suppose $j\sim j^\prime$
in the
sense of Definition 1.3(ii).  Then
$\mj$ is both Dirichlet and Neumann isospectral to $\mjp$.
\endproclaim

The proof is similar to the argument given in \cite{G3}
for the construction of isospectral metrics on nilmanifolds (compact quotients
of
nilpotent Lie groups by discrete subgroups).  Before giving the proof, we
give a geometric interpretation of the condition $j \sim \jp$ and establish
some
notation.  

\definition{1.6 Remarks and Notation}   Suppose $j \sim \jp$.  

(i.)  If $\fz$ is one-dimensional, then $j
\simeq \jp$ in the notation of Definition 1.3, with $C$ being the identity
operator on
$\fz$. Thus the isometry conditions of Proposition 1.4 hold with the isometry
$\tau$ from
$G(j)$ to $G(j')$ given by $\tau(x,z)=(A(x),z)$ with $A$ as in 1.3(i).  If
$\lc$ is any
lattice in
$\fz$, the translations of $G(j)$ and $G(j')$ by elements of $\lc$ commute with
$\tau$, and thus $\tau$ induces global isometries between $\ogj$ and $\ogjp$
and
between $M(j)$ and $M(j')$.

(ii.) If $\fz$ is higher-dimensional, then $\gj$ need not be isometric to
$\gjp$, but
the two manifolds admit many isometric quotients.  More precisely, consider any
co-dimension one subgroup $W$ of the derived group of $\gj$.  Such a subgroup
corresponds under the exponential map to a co-dimension one subspace of $\fz$,
equivalently to the kernel of a non-trivial linear functional $\ll$ on $\fz$. 
Let
$\fz_{\ll}$ be the orthogonal complement of $W$ in $\fz$.  Then the two-step
nilpotent Lie group $G_{\ll}(j) := \gj /W$ with the induced Riemannian metric
is
associated as in 1.2 with the data $(\fv, \fz_{\ll}, j_{|\fz_{\ll}})$.  Observe
that
$j_{|\fz_{\ll}}\sim j^\prime_{|\fz_{\ll}}$ since $j\sim j'$.  Thus  by (i) and
the fact
that
$\fz_{\ll}$ is one-dimensional, we see that 
$\glj$ is isometric to $\gljp$.  

(iii.) If $\lc$ is a lattice in $\fz$ and if $\ll \in \lc^*$, i.e., $\ll$ is
integer-valued on $\lc$, then the projection from $\fz$ to $\fzl$ maps $\lc$ to
a
lattice $\lc_{\ll}$ in $\fzl$.  The associated quotients $\oglj$ and $\ogljp$,
defined
as in 1.2,  are isometric.  Under the identifications in 1.2, the isometry
$\Psi_\lambda$ is given by $\Psi_\ll (x,\bar{\fz})=(A_\ll (x),\bar{\fz})$,
where
$A_\ll \in \so(\fv)$ satisfies $j' (z)=A_\ll j(z)A_\ll^{-1}$ for
$z\in\fz_\ll$.  This isometry restricts to an isometry between the compact
submanifolds
$\mjl$ and $\mjlp$ of $\oglj$ and $\ogljp$ defined as in 1.2(iv).

(iv.)  We will say two vectors $\ll$ and $\mu$ in $\lc^*$ are {\it equivalent},
denoted $\ll \sim \mu$, if they have the same kernel.  Denote the equivalence
class
of $\ll$ by $[\ll]$ and denote the set of equivalence classes by $[\lc^*]$. 
Observe
that $\glj$, $\oglj$ and $\mlj$ depend only on the equivalence class of $\ll$. 
\enddefinition

\proclaim{Lemma 1.7} Let $\pi_\ll :\ogj \to \oglj$ be the homomorphic
projection.  The
Laplacians
$\Delta$ of 
$\ogj$ and $\Delta_\ll$ of $\oglj$ satisfy 
$$\pi_\ll^*\circ\Delta_\ll =
\Delta\circ\pi_\ll^*.$$ 
\endproclaim

\demo{Proof}
It is well-known that the
conclusion holds provided that the
projection is a Riemannian
submersion with totally geodesic
fibers.  The elementary proof that
these conditions hold in our case
is identical to the proof of
Proposition 1.5 in \cite{G3}.
\enddemo 

Note that $\pi_\ll$ gives $\ogj$ the structure of a principal torus bundle. 
Moreover, $\pi_\ll$ restricts to a Riemannian submersion from $\mj$ to $\mlj$
whose
fibers are flat tori.

\remark{1.8 Remark}  $\pi_0$ corresponds to the canonical projection
$\fv\times\bar{\fz}\to \fv$ in the notation of 1.2.  Moreover,
$\overline{G_0(j)}$ is
a flat torus; in fact, it is isometric to the quotient of the Euclidean space
$(\fv ,
<\cdot,\cdot>)$ by a lattice.  The fiber torus is isometric to $(\bar{\fz} ,
<\cdot,\cdot>)$.  In particular the fact that
$\pi_0$ is a Riemannian submersion implies that the Riemannian measure on
$L^2(\ogj)$
coincides with the Lebesgue measure on $\fv\times\bar{\fz}$.
\endremark

\demo{Proof of Theorem 1.5}
In the notation of 1.2, the derived group of the Lie group $\ogj$ is identified
with
the torus $\bar{\fz}$.  This torus
acts isometrically on $\ogj$ and on the
submanifold $\mj$ by left translations.  The
resulting action of $\bar{\fz}$ on $L^2(\mj)$,
given by 
$$(\rho(\bar w)f)(x,\bar z)= f(x,\bar z+\bar
w),\tag1$$ clearly carries the space of smooth
functions with Dirichlet boundary
conditions to itself.  To see that it
also leaves invariant the space of
smooth functions with Neumann boundary
conditions, observe that the normal
derivative of a function $f$ across the boundary of $\mj$ at the point
$(x,\bar{z})$, where $x$ is a
unit vector in $\fv$, is given
by  $xf(x,\bar{z})$ where $xf$ denotes the left-invariant vector field $x$ on
$\ogj$ applied to
$f$. Indeed 
$$xf(x,\bar{z}) =\frac{d}{dt}f((x,\bar{z})\overline{exp}(tx))=
\frac{d}{dt}f((x,\bar{z})(tx,0))=\frac{d}{dt}f((1+t)x,\bar{z})$$
by 1.2(ii),(iii).
Since the torus
$\bar{\fz}$ lies in the center of
$\ogj$, the torus action $\rho$ commutes with all left-invariant vector fields. 
In view of the definition (1) of $\rho$, it follows that $\rho$ leaves
invariant the
space of smooth functions with Neumann boundary conditions.  

By Fourier decomposition on the torus, we can write 
$$L^2(\mj)=L^2(B\times\bar{\fz})=\displaystyle{\sum_{\ll \in\lc^*}\Hl}$$
where
$$\Hl =\{f\in L^2(B\times\bar{\fz}):\rho(\bar z)f=e^{2\pi i\lambda(z)}f \text{
for all
}\bar{z}\in \bar{\fz}\}.$$
By the comments above, the space of smooth functions on $\mj$ with Dirichlet,
respectively Neumann, boundary conditions decomposes into its intersections
with the
$\Hl$.  To avoid cumbersome notation, we will refer to $\spec(\Delta_{|{\Cal
H}_\ll})$ with Dirichlet (or Neumann) boundary conditions to mean the
spectrum of the Laplacian of $M(j)$ restricted to the space of smooth functions
in
$\Hl$ with the appropriate boundary conditions.

Set
$$\hl =\displaystyle{\sum_{\mu\sim\ll}{\Cal H}_\mu}.$$
(See Notation 1.6(iv).)
Define $\Hlp$ and $\hlp$ similarly using the data $(\fv, \fz, j', \lc)$.

By Lemma 1.7 and Remark 1.8, $\pi_0^*$ intertwines the Laplacian $\Delta$ of
$\mj$,
restricted to
${\Cal H}_0$, with the Euclidean Laplacian on the ball $B$ and similarly for
the
Laplacian
$\Delta^\prime$ of
$\mjp$, restricted to ${\Cal H}_0^\prime$.  Thus with either Dirichlet or
Neumann
boundary conditions, we have 
$$\spec(\Delta_{|{\Cal H}_0})=\spec(\Delta^\prime_{|{\Cal
H}_0^\prime}).\tag2$$

Next for $0\neq \ll \in\lc^*$, the map $\pi_\ll^*$ is a unitary map from
$L^2(\mjl)$
to ${\Cal H}_0\oplus\hl$ (i.e., to $\{f\in
L^2(B\times\bar{\fz}):\rho(\bar z)f=f \text{ for all }z\in \ker(\ll)\}$).  Thus
by 1.6(iii) and Lemma 1.7, we have with either Dirichlet or Neumann boundary
conditions
that
$\spec(\Delta_{|{\Cal H}_0\oplus\hl})=\spec(\Delta^\prime_{|{\Cal
H}_0^\prime\oplus\hlp}).$  In view of equation 2, we thus have with either
boundary condition that 
$\spec(\Delta_{|\hl})=\spec(\Delta^\prime_{\hlp})$
for every $\ll\in\lc^*$.  The theorem now
follows.
\enddemo

\remark{1.9 Remarks}(i) The intertwining
operator $T$ between the Laplacians of $\mj$ and
$\mjp$ can be written explicitly as
$T=\oplus_{\ll\in\lc^*}T_\ll$ where $T_\ll
:\Hl\to\Hlp$ is given by $(T_\ll f)(x,\bar{z}) =
f(A_\ll^{-1}(x),\bar{z})$ where $A_\ll$ is given as in 1.6(iii).

(ii) By replacing the ball $B$ with the vector space
$\fv$ everywhere in the argument above, one
obtains a unitary isomorphism
$T:L^2(\ogj)\to L^2(\ogjp)$ satisfying
$\Delta^\prime =T\circ \Delta\circ T^{-1}$,
where $T$ is given by the same formula as
in (i), but with $\Hl$ now being a subspace of
$L^2(\fv\times\bar{\fz})$.

(iii) By working with the Fourier transform on $L^2(\fv\times\fz)$ with respect
to the second variable, one can similarly obtain a unitary
isomorphism between $L^2(\gj)$ and $L^2(\gjp)$
which intertwines the Laplacians.  (There are some technical complications in
the proof; for example, to define $T$, one needs $A_\ll^{-1}(x)$ to be
measurable as a function of $(x,\ll)\in\fv\times\fz^*$.  Note that the $A_z$'s
in Definition 1.3(ii), and thus the $A_\ll$'s in Remark 1.6, are not uniquely
determined.  We have shown that one can choose the $A_z$'s so that the map
$z\to A_z$ from $\fz$ to the orthogonal group $O(\fv)$ is in fact real analytic
on a
Zariski open subset of $\fz$.)  Since
$\gj$ is diffeomorphic to $\rb^n$ for
some $n$, we thus obtain metrics on
$\rb^n$ whose Laplacians are
intertwined.  We omit the details here
as we are currently investigating the
behavior of the scattering operators for
these metrics.  We expect to address
this issue in a later paper.
\endremark

\subhead Example 1.10 \endsubhead  In \cite{G2,3}, examples were given of pairs
of
isospectral (in the sense of Definition 1.3), inequivalent linear
maps $j, j^\prime : \fz\ra \so({\frak v})$, where $\fz$ was
$3$-dimensional and ${\frak v}$ was $4n$-dimensional with $n\geq 2$.  The
resulting isospectral manifolds, given by Theorem 1.5, thus have
minimum dimension eleven.  (The fact that $j$ and $j^\prime$ give rise to
isospectral compact manifolds with boundary was not observed in \cite{G2,3}. 
Instead $j$ and $j^\prime$ were used to construct isospectral closed manifolds
using
the method described in \S 3 below.)  We now construct
$7$-dimensional examples.  As we'll see in \S 4, these have quite different
geometric
properties from the earlier examples. 

Let $H$ be the quaternions and $P$ the pure quaternions, i.e.,
$P=\{q\in H:\bar{q}=-q\}$.  For $q\in H$, let $L(q)$ and $R(q)$
denote left and right multiplication by $q$ on $H$.  For $q,p\in
P$, set $J(q,p)=L(q) +R(p)$.  Then $J(q,p)$ is skew-symmetric
relative to the standard inner product on $H$.  Indeed the
decomposition $\so(4)=\so(3) + \so(3)$ says that all skew-symmetric
operators are of this form.  An easy computation shows that the
eigenvalues of $J(q,p)$ are $\pm i\sqrt{|q|^2+|p|^2}$ and $\pm
i\sqrt{|(|q|^2-|p|^2)|}$; in particular, the spectrum of $J(q,p)$
depends only on the lengths of $q$ and $p$.

Now let ${\frak v}=H$, viewed as $\rb^4$ with the standard inner product, and
let $\fz=P$, viewed as $\rb^3$ with the standard inner product.  Let
$T$ and $\tp$ be fixed invertible linear operators on $P$ such
that $\tp=A\circ T$ where $A$ is an orthogonal operator of
determinant $-1$.  Define $j,j^\prime :\fz\ra\so({\frak v})$ by
$j(q)=J(q,Tq)$ and $j^\prime(q)=J(q,T^\prime q)$ for all $q\in P$. 
Then $j\sim j^\prime$ in the sense of Definition 1.3.

We next check whether $j$ is equivalent to $j^\prime$.  The group
SO$({\frak v})$ consists of all operators $L(a) R(b)$ where $a$ and $b$ are
unit quaternions.  Conjugation of $J(q,p)$ by $L_a R_b$ gives
$J(a^{-1} q a, b p b^{-1})$.  All orthogonal transformations of ${\frak v}$
are compositions of elements of SO$({\frak v})$ with the quaternionic
conjugation map $B$ of ${\frak v}$.  Conjugation of $J(q,p)$
by $B$ yields $J(-p,-q)$.  Since $\det\,T^\prime=-\det\,(T)$, it follows
easily that the construction above always yields inequivalent maps
$j$ and $j^\prime$.

With any choice of lattice ${\Cal L}$ in $\fz$, Theorem 1.5
yields pairs of isospectral 7-dimensional compact manifolds with boundary which
are
not locally isometric.

\head\S2 Examples of Isospectral Lie Algebra Deformations\endhead

\definition{Definition 2.1}
Let ${\frak v}$ and $\fz$ be finite dimensional inner product spaces and $j_0$
any linear map from $\fz$ to $\so({\frak v})$.  By a $d$-{\it parameter
non-trivial
isospectral deformation} of $j_0$ we mean a continuous function $u\mapsto\ju$
from a
pathwise connected subset $D$ of $\rb^{d}$ having non-empty interior into the
space of linear maps from $\fz$ to $\so({\frak v})$ such that
\smallskip\roster
\item"(i) " $\jo=j_{u_{0}}$ for some $u_0\in D$;
\item"(ii) " $j_u\sim j_0$ for all $u\in D$ (see Definition
1.3(ii));
\item"(iii) " $j_{u}\not\simeq j_{u^\prime}$ whenever $u$ and
$u^\prime$ are distinct points in $D$ (see Definition 1.3(i)).
\endroster
\enddefinition

Equivalently, ${\Cal G}=
\{{\frak g}(j_u):u\in D\}$ is a family containing ${\frak g}(j_0)$ of nilpotent
metric
Lie algebras all
having ${\frak v}\oplus\fz$ as their underlying vector space, and the structure
constants of ${\frak g}(j_{u})$ relative to any fixed bases of ${\frak v}$ and
$\fz$
depend continuously on the parameter $d$-tuple $u$.  Any choice of lattice
${\Cal L}$
of maximal rank in $\fz$ gives rise to a $d$-parameter family $\{M(j_u)\}_{u\in
D}$ of
isospectral compact manifolds with boundary as in Theorem 1.5.

Throughout this section, we will consider the special case where
$\dim\fz=2$ with $m=\dim {\frak v}$ variable.  Our goal is to show that when
either 
$m=5$  or $m \geq 7$, every ``generic'' ${\frak j}_0$ admits a
$d$-parameter non-trivial isospectral deformation with $d> 1$.  For
$m=6$, we will exhibit explicitly one parameter deformations for
certain $j_0$ of a restrictive type.  For $m\leq 4$,
straightforward calculations show that any two isospectral $j$'s
are in fact equivalent, so non-trivial isospectral deformations of this type
are
impossible.

\proclaim{Theorem 2.2}Let $\dim\fz=2$, let $m=\dim {\frak v}$ be any positive
integer
other than
$1,2,3,4$, or $6$, and let $L$ be the real vector space consisting of all
linear maps
from $\fz$ to $\so({\frak v})$.  Then there is a Zariski open subset ${\Cal O}$
of $L$ (i.e., ${\Cal O}$ is the complement of the set of roots of some
non-zero polynomial function on $L$) such that each $\jo\in {\Cal
O}$ admits a $d$-parameter non-trivial
isospectral deformation where $d\geq m(m-1)/2 - [m/2]([m/2]+2)>1$.
\endproclaim

\demo{Proof} 
For $j\in L$, let $$I_j=\{\jp\in L:\jp \sim j\}
\text{ and }E_j=\{\jp\in I_j:\jp \simeq j \}.\tag1$$  The idea of the
proof is to define ${\Cal O}$ in such a way that for
$\jo\in\oc$, $P_{j_{0}}:=I_{j_{0}}\cap\oc$ is an embedded submanifold
of $L$ which can be foliated by its intersection with the sets
$E_j, j\in P_{j_{0}}$, and for which there is a submanifold
$N_{j_{0}}$ of $P_{j_{0}}$ transverse to the foliation.  Any parametrization of
$N_{j_{0}}$ then defines a non-trivial isospectral deformation of $j_0$.

Let $l=[m/2]$ and, for $1\leq k\leq l$, define $T_k:\so(\fv)\ra\rb$
by $T_k(C)=\trace(C^{2k})$.  If $C$ and $C^\prime$ are
similar, i.e., have the same eigenvalues, trivially
$T_k(C)=T_k(C^\prime)$ for all $k$.  But the converse is also true as
can be seen by a standard combinatoric argument showing that the
coefficients of
powers of $\la$ in the characteristic polynomial
$\cx(\la,C)=\det(\la Id-C)$ are polynomials in
$\{T_1(C),\ldots,T_\ell(C)\}$.  If we define $T_k:\fz\times L\ra \rb$
by $T_k(z,j)=T_k(j(z))$, this means that
$j\sim \jp \Leftrightarrow T_k(z,j)=T_k(z,\jp)$ for
all $z\in\fz$ and all $k, 1\leq k\leq\ell$.  Moreover, each of the
functions $T_k$ is a polynomial on $\fz\times L$ which is separately
homogeneous of degree $2k$ in each variable.  If we fix any
orthornormal basis $\{\epsilon_1,\epsilon_2\}$ of $\fz$ and denote a typical
element $z\in\fz$ by $z=s\epsilon_1+t\epsilon_2$, expansion of
$T_k(z,j)={\trace}(sj(\ei)+tj(\ez))^{2k}$ into
$(s,t)$
monomials gives us $2{k+1}$ coefficient functions which are
polynomials in $j(\epsilon_1)$ and $j(\epsilon_2)$ and thus polynomials on
$L$.  Since $\sum_{k=1}^{\ell}(2k+1)=\ell(\ell + 2)$, we conclude
that there is a map $F:L\ra\rb^{\ell(\ell+2)}$ each of whose entries
is a polynomial on $L$ and for which
$j \sim \jp \Leftrightarrow F(j)=F(\jp)$.  Let $R$ be the
maximum rank of $F$.  Thus $R$ is the largest integer
for which there is some $j\in L$ such that the tangent map
$F_{*_{j}}:L\ra\rb^{\ell(\ell +2)}$ has rank $R$.  Since each of the
entries in any matrix representation of $F_{*_{j}}$ is a polynomial in
$j$ and since a matrix has rank $\geq R$ precisely when the sum of
the squares of the determinants of its $R\times R$ minors is non-zero,
it follows that the subset ${\Cal O}_1$ of $L$ on which $F$ has rank
$R$ is a Zariski open set.  Moreover, for $\jo\in\oc_{1}$, the Implicit
Function
Theorem says that there is a neighborhood
${\Cal U}$ of $\jo$ in $L$ for which $I_{\jo}\cap
U=F^{-1}(F(\jo))\cap U$ is an embedded submanifold
of co-dimension $R$.

We now turn toward examination of the sets $E_j$ in (1).  The group
$G=O(\fv)\times O(\fz)$ acts on $L$ by $$((A,C)\cdot j)(z)=A
j(C^{-1}z)A^{-1}\tag2$$ and, by Definition 1.3, $\jp \simeq
j \Leftrightarrow\jp=(A,C)\cdot j$ for some $(A,C)\in G$.  Let
$1_{\fv}$ and $1_\fz$ denote the identity operators on $\fv$ and $\fz$.
 We now claim that there is a Zariski open subset $\oc_2$ of $L$
such that for each $j\in\oc_2, E_j$ is the orbit of $j$ under the
subgroup $K:=O(\fv)\times \{\pm 1_\fz\}$ and such that the stability subgroup
of
$K$ at $j$ is $\{(\pm 1_\fv, 1_\fz)\}$.  To see this, first consider
any $j\in L$ and $(A,C)\in G$ such that $(A,C)\cdot j\in E_j$.  Since
$E_j\subset
I_j$ and since
$j\circ C^{-1}=(A^{-1},1)\cdot(A,C)
\cdot j$, we see that $j\circ C^{-1}\in I_j$. Thus $T_k(z,j\circ
C^{-1})=T_k(z,j)$ for
all
$z$ and
$k$.  In particular, $C$ is orthogonal both with respect to the given inner
product on $\fz$ and the quadratic form
$z\mapsto T_1(z,j)=\trace(j(z))^2= -<j(z),j(z)>$, where
$<c,d>=\trace(c \,d^t)= -\trace(c d)$ is the standard inner product on
$\so(\fv)$.  Relative to any orthonormal basis
$\{\epsilon_1,\epsilon_2\}$ of $\fz, T_1(\cdot,j)$ has matrix 
$$-\left[\matrix |j(\epsilon_1)|^2 &<j(\epsilon_{1}),j(\epsilon_2)>
\\                         <j(\epsilon_1),j(\epsilon_2)>&|j(\epsilon_2)|^2
\endmatrix\right].
$$
Unless this matrix is a scalar multiple of the $2\times 2$ identity
matrix, i.e. unless $$\phi_1(j):=(|(j(\ei)|^{2}-|j(\e2)|^2)^2 +
<j(\ei),j(\e2)>^2$$ vanishes, there are precisely four transformations
orthogonal with respect to both forms, namely $\pm 1_{\fz}$ and $\pm C_0$
where $C_0$ is the reflection leaving one eigenvector of the above matrix
fixed while changing the sign of the other.  For $j$ not a root of the
polynomial $\phi_1$, we conclude that $C$ is one of these four
transformations.  But $m\geq 5$ means $\ell=[m/2]\geq 2$ so $C$ must also
satisfy $T_2(z,j\circ C^{-1})=T_2(z,j)$, i.e.
$\trace(j(z))^4=\trace(j(C^{-1} z))^4$ for all $z\in\fz$.  By
straightforward but tedious calculations, one can check that there is a
fourth order polynomial $\phi_2$ on $L$ for which
$T_2(\cdot,j\circ C_0^{-1})\neq T_2(\cdot, j)$ when
$\phi_2(j)\neq0$.  Thus when both $\phi_1(j)$ and $\phi_2(j)$
are non-zero, $(A,C)\cdot j\in E_j\Leftrightarrow (A,C)=(A,\pm 1_\fz)\in
K$.  In this case, $(A,C)\cdot j=j$ if and only if either $C=1_\fz$ and $A\in
O(\fv)$
commutes with $j(z)$ for all $z$, or else $C=-1_\fz$ and $A$ anti-commutes with
$j(z)$ for all $z$.  With $\{\ei,\e2\}$ as above and $j_1=j(\ei),
j_2=j(\e2)$, it's easy to select choices of $j_1$ and $j_2$ for which no
non-zero linear operator $A$ on $\fv$ anti-commutes with both $j_1$ and
$j_2$, while $\pm 1_\fv$ are the only orthogonal operators commuting with
both $j_1$ and $j_2$.  Moreover, these properties are equivalent to saying
that the linear map $\phi_j(B):=(j_1 B-Bj_1, j_2 B-Bj_2)$ from $\fgl(\fv)$
to $\fgl(\fv)\times\fgl(\fv)$ has one-dimensional kernel
while $\widetilde{\phi}_j(B):=(j_1 B +Bj_1, j_2 B+Bj_2)$ is injective,
conditions which can be expressed by the statement that certain non-vanishing
polynomials $\phi$ and
$\widetilde\phi$ on $L$ do not have $j$ as a root.  Combining all of these
arguments, when $j$ belongs to
the complement $\oc_2$ of the set of roots of
$\phi_1^2+\phi_2^2+\phi^2+\widetilde{\phi}^2$, the
properties announced above in our claim are satisfied.

Let $\oc=\oc_1\cap\oc_2, \jo\in\oc$, and $P_{\jo}=I_{\jo}\wedge\oc$.  From 
above, $P_{\jo}$ is a smooth manifold whose dimension is $\dim
L-R\geq m(m-1)-[m/2]([m/2] +2)$.  For $K=O(m)\times\{\pm 1_\fz\}$,
it's trivial to check that when any one of the polynomials defining
${\Cal O}$ does not vanish at $j$, the same is true for each member of
the $m(m-1)/2$-dimensional orbit $K\cdot j$; i.e., ${\Cal O}$ is closed under
the
action of $K$.  Moreover, for
$j\in\oc$, we have shown that the orbit $K\cdot j$ coincides with $E_j$ and the
stability subgroup at $j$ is
$Z:=\{(\pm 1_{\fv}, 1_{\fz})\}$.  This means that the compact group
$K/Z$ acts freely on the manifold $P_{\jo}$ with orbits expressing
equivalence of elements.  By the properties of compact transformation
groups (e.g. \cite{Bd}, pp 82--86), there is a submanifold $N_{\jo}$ of
$P_{\jo}$ such that $(j,\tilde{K})\mapsto \tilde{K}(j)$ is a
homeomorphism from $N_{\jo}\times (K/Z)$ onto an open neighborhood of
$\jo$ in $P_{\jo}$.  The dimension of $N_{\jo}$ is then
$$
d=\dim I_{\jo}-m(m-1)/2\geq \frac{m(m-1)}{2}-[m/2]([m/2]+2).
$$  For $m=5$ or $m\geq 7$, clearly $d> 1$ and any local parameterization of
$N_{\jo}$ defines a $d$-parameter non-trivial isospectral deformation
of $\jo$.  
\enddemo

\bigskip
\subhead 2.3 Eight dimensional examples\endsubhead    For $m=6$, the
argument in the proof of Theorem 2.2 breaks down since
$m(m-1)/2-[m/2]([m/2]+2)=15-15=0$.  In the language of the proof of Theorem
2.2, the
examples below correspond to choosing certain $\jo$'s where the rank of the
polynomial
map $F$ is less than $R$ with the result being that the isospectral surface
$I_{\jo}$ in equation (1) is four-dimensional while the sets $E_j$ contained in
$I_{\jo}$ are three-dimensional and admit a one-parameter
transversal.  Lengthy and non-illuminating calculations are avoided
by fixing orthonormal bases for $\fv$ and $\fz$ and simply defining in
concrete matrix terms the members of the transversal.

Thus, take $\fz=\rb^2$ and $ \fv=\rb^6$ with their standard ordered bases
and standard inner product.  For $a,b\in\so(6)$ and $s,t\in\rb$, define
$j_{a,b}(s,t)=sa+tb$.  Each linear map
$j:\rb^2\ra\so(6)$ is of the form $j=j_{a,b}$ for some $a,b\in \so(6)$. 
Fix for the remainder of the discussion an element $a\in\so(6)$ which is in
block
diagonal form with $2\times 2$ diagonal blocks $a_i\left[\matrix 0&-1\\
1&0\endmatrix\right], 1 \leq i
\leq 3$, where $0<a_1<a_2<a_3$.  Consider all matrices $b\in \so(6)$ of the
form
$$ b=\left[\matrix 0&0&b_{12}&0&b_{13}&0\\
                             0&0&0&0&0&0\\
                             -b_{12}&0&0&0&b_{23}&0\\
                             0&0&0&0&0&0\\
                             -b_{13}&0&-b_{23}&0&0&0\\
                             0&0&0&0&0&0\endmatrix\right]
$$
with $(b_{12},b_{13},b_{23}) \in\rb^3-\{0\}$.

We first note that if $b$ and $b^\prime$ are of this form, then
$j_{a,b}\simeq j_{a,b^\prime}\Leftrightarrow b^\prime=\pm b$. 
Indeed, in the notation of equation (2), if $j_{a,\bp}=(A,C) \cdot j_{a,b}$ for
some
$A\in O(6), C\in O(2)$, then for $\e2=(0,1), j_{a,b^\prime}(\e2)=b^\prime$ is a
rank
$2$ matrix similar to $j_{a,b}(C^{-1} \e2)$.  But a simple calculation shows
that
$j_{a,b}(s,t)$ has rank 2 only when $s=0$.  It follows that
$C\e2=\pm\e2$, so $C$ is one of $\pm\left[\matrix 1&0\\ 0&1\endmatrix\right]$,
$
\pm\left[\matrix 1&0\\ 0&-1\endmatrix\right]$, and then $Aa A^{-1}=\pm a$, $A b 
A^{-1}=\pm\bp$.  Since $a_1, a_2, a_3$ are distinct, this forces $A$ to be in
block
diagonal form with $2\times 2$ diagonal blocks which either all commute with
$\left[\matrix 0&1\\-1&0\endmatrix\right]$ or all anticommute with
$\left[\matrix 0&1\\-1&0\endmatrix\right]$.
Using the specific form of $b$ and $b^\prime$, it follows in either case
that $A b A^{-1}=b$ so $b^\prime=\pm b$.

Next an easy direct calculation yields
$$\det\{\la Id -j_{a,b}(s,t)\}=\prod_{i=1}^{3}(\la^2 + s^2 a_i^2)+\la^4
t^2\sum_{i<j} b_{ij}^2+\la^2 s^2 t^2(a_1^2 b_{23}^2 +a_2^2 b_{13}^2
+a_3^2 b_{12}^2).$$  Comparing coefficients, it follows that
$j_{a,b}\sim j_{a,\bp}\Leftrightarrow (b_{12}, b_{13}, b_{23})$
and $(b_{12}^\prime, b_{13}^\prime, b_{23}^\prime)$ satisfy the equations
$$
\sum_{i<j} b_{ij}^2-(b_{ij}^\prime)^2=0\tag{i}
$$
and
$$a_1^2 (b_{23}^2-(b_{23}')^2) +a_2^2 (b_{13}^2-(b_{13}')^2)
+a_3^2(b_{12}^2-(b_{12}')^2)=0.\tag{ii}
$$
In view of equation (i), equation (ii) can be rewritten as
$$\sum_{i<j}(a_i^2+a_j^2)(b_{ij}^2-(b_{ij}^\prime)^2)=0.\tag{ii'}
$$
The general solution of equations (i) and (ii') is
$$\align (b_{12}^\prime)^2&=b_{12}^2 +u(a_2^2-a_1^2)\\
(b_{13}^\prime)^2&=b_{13}^2
+u(a_1^2-a_3^2)\tag*\\ (b_{23}^\prime)^2&=b_{23}^2+u(a_3^2-a_2^2)\endalign$$
where $u$ is any real number in the closed interval
$I=[\max\left(\frac{-b_{12}^2}{a_2^2-a_1^2},\frac{-b_{23}^2}{a_3^2-a_2^2}
\right),\frac{b_{13}^2}{a_3^2-a_1^{2}}]$.  If we take any $b$ for
which $I$ has non-empty interior and, for each $u\in I$, define
$b(u)$ as the unique solution of the above equations for which
$b_{ij}(u)$ has the same sign as $b_{ij}$ for all $i,j$, it follows
that $u\ra j_{a,b(u)}$ is a $1$-parameter non-trivial isospectral
deformation of $j_{a,b}$.
\head \S3 Compact Nilmanifolds\endhead  

A compact Riemannian nilmanifold is a quotient
$N = \Gamma \backslash G$ of a simply-connected nilpotent Lie group $G$ 
by a (possibly trivial) discrete subgroup $\Gamma$, together with a Riemannian
metric $g$ whose lift to $G$ is left-invariant.

We now recall a method, developed in \cite{G3}, for constructing isospectral
compact
Riemannian nilmanifolds.  For convenience, we'll restrict our attention to
two-step
nilmanifolds, although Theorem 3.2 below can be formulated in the context of
nilmanifolds of arbitrary step.   Even in the two-step case, the formulation of
Theorem 3.2 in
\cite{G3, Theorem 1.8} is slightly more general than that given here.   

\definition{Notation and Remarks 3.1} (i) A nilpotent Lie group $G$ admits
a co-compact discrete subgroup $\Gamma$ if and only if the Lie algebra ${\frak
g}$ of
$G$ has a basis ${\Cal B}$ relative to which the constants of structure are
integers
(see
\cite{R}). If ${\Cal B}$ is such a basis and ${\Cal A}$ is the integer span of
${\Cal B}$, then $\exp({\Cal A})$ generates a co-compact discrete subgroup of
$G$. 
Conversely, if $\Gamma$ is a co-compact discrete subgroup of $G$, then
$\log(\Gamma)$
spans a lattice of full rank in $\fg$, where $\log:G\to\fg$ is the inverse of
the Lie
group exponential map.
  
(ii) We use the notation of 1.1 and 1.2.  Thus a simply-connected nilpotent Lie
group 
$G=G(j)$ with a left-invariant metric is defined by  data $(\fv, \fz\, j)$.  If
$\Gamma$ is a co-compact discrete subgroup of $G$, then $\Gamma$ intersects
$[G,G]$
in a lattice of full rank ${\Cal L}$, which we may also view as a lattice in
$\fz$
under the identification in 1.2. In
summary, a compact nilmanifold
$N=\Gamma\backslash G$ is defined by the data $(\fv, \fz,j,\Gamma)$ and
$\Gamma$
determines a lattice $\lc$ in
$\fz$.  In the sequel, we will consider fixed $(\fv,\fz,\lc)$ but vary the
choice of
$j$ with the requirement that the resulting simply-connected nilpotent Lie
group
$G(j)$ admit a co-compact discrete subgroup $\Gamma$ whose intersection with
the
derived group of $G(j)$ is given by $\lc$.  We will denote the nilmanifold
$\Gamma\backslash G(j)$ by $\njg$.
 
(iii) We continue to use the notation of 1.6 as well.  For $\ll\in\lc^*$, the
projection $G(j)\to G_\ll(j) =G(j)/\ker(\ll)$ sends $\Gamma$ to a co-compact
discrete
subgroup $\Gamma_\ll$.  We denote by $\nlj$ the quotient
$\Gamma_\ll\backslash G_\ll(j)$ with the Riemannian metric induced by that of
$G_\ll(j)$.  Note that $N_0(j,\Gamma)$ is a flat torus.  Letting ${\Cal A}_\fv$
be
the image of $\log(\Gamma)$ under the orthogonal projection from $\fg(j)$ to
$\fv$,
then
$N_0(j,\Gamma)$ is isometric to the torus $\fv/{\Cal A}_\fv$ with the flat
metric
defined by the inner product on $\fv$.
\enddefinition

\proclaim{Theorem 3.2 \cite{G3}}  Let $\njg$ and
$\njp$ be compact Riemannian
nilmanifolds associated with the data $(\fv, \fz,\lc)$ as in 3.1. Suppose that
$\spec(\nlj) =
\spec(\nljp)$ for every
$\lambda \in \lc^{*}$. Then $\spec(\njg) = \spec(\njp)$.
\endproclaim

(We wish to correct an error in the version of this theorem given in \cite{G3},
Theorem
1.8:  One must assume that the correspondence $\ll\to\ll^\prime$ given there is
norm-preserving.  This assumption is actually satisfied in all the applications
of
Theorem 1.8 given in \cite{G3}.)

\definition{Definition 3.3} We will say a two-step nilpotent Lie group $G=G(j)$
is
{\it non-singular} if
$j(z)$ is non-singular for all $z \in \fz$. We will also say any associated
compact
nilmanifold $N(j,\Gamma)$ is non-singular in this
case.
\enddefinition

In \cite{G3}, we studied non-singular nilmanifolds and proved the following as
a
consequence of Theorem 3.2:

\proclaim{Theorem 3.4}  In the notation of 3.1, let $\njg$ and $\njp$
be compact non-singular two-step Riemannian nilmanifolds associated with the
same data
$(\fv,\fz,\lc)$. Assume:

\item{(i)} $\spec(N_{0}(j,\Gamma)) = \spec(N_{0}(j^{\prime},\Gamma^\prime))$
and
\item{(ii)} $j\sim j^\prime$. (See Definition 1.3(ii).)

\noindent Then $\spec(\njg)=
spec(\njp)$.
\endproclaim

\subhead Example 3.5 \endsubhead Examples of minimum dimension $11$ were given
in
\cite{G2},
\cite{G3}.  We now construct compact quotients
of the pairs of 7-dimensional isospectral simply-connected manifolds $G(j)$ and
$G(j^{\prime})$ constructed in Example 1.10.

In the notation of Example 1.10, observe that the constants of 
structure of ${\frak g}(j)$ relative to the ``standard'' basis are integers
provided
that the matrix entries of $T$ relative to the standard basis of ${\frak z}$
are integers. Thus we assume that the matrix entries of both $T$ and
$T^{\prime}$ are
integers. We can then, for example, let ${\Cal A}$ and ${\Cal A}^\prime$ be the
integer
span of the  standard basis elements of ${\frak v}$ and ${\frak z}$ and let
$\Gamma$
and $\Gamma '$ be the discrete subgroups of $G(j)$ and $G(j')$ generated by
$\exp({\Cal A})$ and $\exp '({\Cal A}')$, respectively.  (See 3.1(i).)  The 
nilmanifolds
$\njg$ and
$\njp$ trivially satisfy condition (i) of Theorem 3.4.; in fact, the tori
$N_0(j,\Gamma)$ and
$N_0(j',\Gamma ')$ are isometric. Moreover, condition (ii) of Theorem 3.4 is
automatic
from the construction in 1.10.  Thus the nilmanifolds
$\njg$ and
$\njp$ are isospectral.
\medskip

The non-singular compact nilmanifolds $N$ are particularly easy to work with as
the
quotient manifolds $\nlj$ defined in 3.1(iii) are Heisenberg manifolds
when $\lambda \neq 0$; that is, the center of $G_{\lambda}$ is one-dimensional.
In \cite{GW2}, the authors gave sufficient conditions for two Heisenberg
manifolds
to be isospectral. (Pesce \cite{P2} later proved these conditions are also
necessary.)
These conditions are used in \cite{G3} to prove Theorem 3.4.

We want to find isospectral compact quotients of some pairs of simply-connected
nilpotent Lie groups associated with the Lie algebras constructed in \S 2. Thus
we need
to generalize Theorem 3.4 to the possibly singular case. As always, we will
assume that
$j(z)$ is non-zero for all $z
\in{\frak z}$ (i.e., that ${\frak z} = [{\frak g}(j), {\frak g}(j)])$. When
$\la\neq
0$, the quotient
${\frak g}_{\lambda}(j)$, has one-dimensional derived algebra but may
have a higher-dimensional center. The corresponding Lie group is of the form 
$G_\la(j) = H \times A$, where $H$ is a Heisenberg group and $A$ an abelian
group. 
Thus in view of Theorem 3.2, we first  need  to examine isospectrality
conditions
for compact quotients of groups of this form.

\definition{Notation 3.6}  In the notation of 3.1, consider a nilmanifold
$\njg$ with
$\fz$ one-dimen\-sional.  We can write
$\fv$ as an orthogonal direct sum
$\fv ={\frak u}\oplus{\frak a}$ where ${\frak a}=\ker(j(z))$ for $0\neq z\in
\fz$. 
(Note that
${\frak a}$ is independent of the choice of $z$ since $\fz$ is
one-dimensional.)  The
Lie algebra
$\fg(j)$ then splits into an orthogonal sum of ideals ${\frak h}\oplus {\frak
a}$,
where ${\frak h}={\frak u}+\fz$ is a Heisenberg algebra.

Since
${\frak a} + {\frak z}$
is the center of ${\frak g}(j)$, $\log(\Gamma)$ intersects 
${\frak a} + {\frak z}$ in a lattice
${\Cal K}$ of maximal rank.  (See \cite{R}.)
Let ${\Cal K}^*$ denote the dual lattice in $({\frak a} + {\frak z})^*$.  The
inner
product
$<\,,\, >$ on ${\frak a} + {\frak z}$ defines a dual inner product
on $({\frak a} + {\frak z})^*$
and thus defines a norm $\| \ \|$ on ${\Cal K}^*$.   
\enddefinition

\proclaim{Proposition 3.7}  Using the notation of 3.1 and 3.6, let $\njg$ be a
compact nilmanifold and assume $\fz$ is one-dimensional. Then
$\spec(\njg)$ is completely determined by the following data:  

\item{(i)} $\spec(N_{0}(j,\Gamma))$.
\item{(ii)} the eigenvalues of the linear operator $j(z)$, where $z$ is a unit
vector in ${\frak z}$. (Since ${\frak z}$ is one-dimensional and $j(z)$ is
skew,
the eigenvalues of $j(z)$ are independent of the choice of unit vector
$z$.)
\item{(iii)} $\{(\sigma (z), \| \sigma \|) \in \rb^{2} : \sigma \in
{\Cal K}^*\}$
where $z$ is given as in (ii).

\endproclaim

The case in which $M$ is a Heisenberg manifold is proven in
\cite{GW2} and is the key lemma used in Theorem 3.4 above.  Proposition 3.7
will be
proved in the Appendix.

\remark{3.8 Remark} In the special case that ${\Cal K}=({\Cal K} \cap {\frak
z})
\oplus ({\Cal K} \cap {\frak a})$,
the data (iii) can be expressed more simply. The inner product $<\,,\,>$
defines flat
Riemannian metrics on the circle ${\frak z} /({\Cal K} \cap {\frak z})$ and the
torus
${\frak a}/({\frak a} \cap {\Cal K})$.  Specifying the data (iii) is equivalent
to 
specifying the length of this circle and the spectrum of this torus. 

\endremark

\medskip

\definition{Notation and remarks 3.9} In 2.3, we considered a class of
eight-dimensional  metric Lie algebras ${\frak g}(j_{a,b})$.
We now show that for certain choices of pairs $j=j_{a,b}$ and 
$j^{\prime}=j_{a,b^\prime}$, the associated
nilpotent Lie groups $G(j)$ and $G(\jp)$ admit isospectral compact quotients. 
 First observe that if the matrix
entries $a_1,a_2,a_3$ of $a$ and
$b_{12},b_{13}, b_{23}$ of $b$ are integers, then the constants of structure of
${\frak g}(j_{a,b})$ with respect
to the standard bases
$\{e_1,\ldots,e_6\}$ of $\fv=\rb^6$ and
$\{\epsilon_1,\epsilon_2\}$ of $\fz=\rb^2$ are integers.  Thus, if we let
${\Cal A}$ be
the lattice in $\fv +\fz$ spanned by
$\{e_1,\ldots,e_6, \epsilon_1,\epsilon_2\}$, then $\exp({\Cal A})$ generates a
co-compact discrete subgroup
$\Gamma_{a,b}$ of
$G(j_{a,b})$.  (See 3.1)
\enddefinition

\proclaim{Theorem 3.10}  In the notation of 2.3 and 3.9, assume that the matrix
entries
of
$a,b$, and $b^\prime$ are integers, that
$\g.c.d.(b_{12}, b_{13}, b_{23})=\g.c.d.(b_{12}^\prime, b_{13}^\prime,
b_{23}^\prime)$, and that the condition (*) of 2.3 is satisfied.  Then the
compact Riemannian nilmanifolds $N(j_{a,b},\Gamma_{a,b})$ and
$N(j_{a,b'},\Gamma_{a,b'})$ are isospectral.
\endproclaim

\demo{Proof}  Write $N=N(j_{a,b},\Gamma_{a,b})$ and
$N'=N(j_{a,b'},\Gamma_{a,b'})$.
 We apply Theorem 3.2.  In the notation of 3.1(ii), the lattice in
$\fz=\rb^2$ associated with both $\Gamma_{a,b}$ and $\Gamma_{a,b'}$ is given by
$\lc
=\span_{\bold Z}
\{\epsilon_1,\epsilon_2\}={\bold Z}^2$.

Let $\la\in \lc^*$.  In case $\ll=0$, both $N_0$ and
$N'_0$ are isometric to the $6$-dimensional cubical torus, so
$\spec (N_0)=\spec(N'_0)$ holds trivially.  Next, observe that
$j(s\epsilon_1+t\epsilon_2)=sa+tb$ is non-singular except when $s=0$.  Thus, if
$\la(\epsilon_1)\neq 0$, then $j_{a,b}(z)$ and $j_{a,b'}(z)$ are non-singular
similar operators for all $z\in\fz_\ll$,  and therefore $N_\la$ and
$N'_\la$ are Heisenberg manifolds.  Proposition 3.7 (see the simplified version
3.8 with
${\frak a}=0$) implies
$\spec\,(N_\la)=\spec\,(N'_{\la})$.

It remains to consider the case $\ker (\la)=\rb\epsilon_1$.  In this case,
$G_\la(j_{a,b})$ and
$G_\la(j_{a,b'})$ are isomorphic as Lie groups to $H\times A$, where $H$ is the
$3$-dimensional Heisenberg group and
$A=\rb^4$.  Letting $\pi_\la:G(j_{a,b})\ra G_\la(j_{a,b})$ be the projection
and
writing
$\overline{X}=\pi_{\la *}(X)$ for
$X$ in the Lie algebra ${\frak g}(j_{a,b})$, we have in the notation of 1.6,
3.6, and
3.9 that
$$
\align {\frak a}&=\pa(\ker (b))=\span\{\overline{e}_2,
\overline{e}_4,\overline{e}_6, b_{23}\overline{e}_1 -b_{13}\overline{e}_3
+b_{12}\overline{e}_5\},\\ {\frak u}&=\pa(\fv\ominus \ker (b)),\endalign
$$ and $\fz_\la=\rb\overline{\epsilon}_2$.  Moreover, letting
$\kc=\pi_{\la *}({\Cal A})\cap({\frak a}+\fz_\la)$, we have
$$
\kc=(\kc\cap \fz_\la)\oplus(\kc\cap{\frak a})
$$ with $\kc\cap\fz_\la={\bold Z}\overline{\epsilon}_2$ and
$\kc\cap {\frak a}=\span_{{\bold Z}}\{\overline{e}_2,
\overline{e}_4,\overline{e}_6,\overline{w}\}$ where 
$$\overline{w}=\g.c.d.(b_{12},b_{13},b_{23})(b_{23}\overline{e}_1 - b_{13}
\overline{e}_3
+ b_{12}\overline{e}_5).
$$ Thus $\kc$ is an orthogonal lattice isomorphic to ${\bold Z}^4\times
|\overline{w}|{\bold Z}$.

The analogous statements hold of course when $b$ is replaced by $b'$.  The data
(ii) in
Proposition 3.7 agree for
$N_\la$ and $N'_\la$; both are given by the eigenvalues of $b$. (Recall that
$b$ and
$b'$ are similar.)  To see that the data (iii), as simplified in Remark 3.8,
agree for
$N_\la$ and
$N'_\la$, we need only show that
$|\overline{w}|=|\overline{w}^\prime|$.  This equality follows from the
hypothesis of the theorem and the fact
that $b_{12}^2 +b_{13}^2 +b_{23}^2=(b_{12}^\prime)^2+(b_{13}^\prime)^2
+(b_{23}^\prime)^2$, as can be seen from
the isospectrality condition $(*)$ of 2.3.
\enddemo

\subhead Example 3.11\endsubhead Fix a choice of $a$ with integer entries
$a_1,a_2,a_3$.  It is easy to find pairs $b$ and $b^\prime$ with integer
entries $b_{ij}$ and $b_{ij}^\prime,
1\leq i < j\leq 3$, so that the isospectrality condition $(*)$ in 2.3 holds,
i.e., so that
$j_{a,b}\sim j_{a,b^\prime}$.  We need only choose the parameter $u$ in $(*)$
so that each of $u(a_i^2 -a_j^2)$
is a difference of two squares; i.e, each $u(a_i^2-a_j^2)$ is an integer
congruent to $0,1$, or $3$ mod $4$.  For
a specific example, take $a_1=1, a_2=2, a_3=3$ and
$u=3$.  We can then take $b_{12}=4, b_{13}=7, b_{23}=7, b_{12}^\prime=5,
b_{13}^\prime=5, b_{23}^\prime=8$. 
In this example, the hypothesis of Theorem 3.10 concerning the g.c.d. of the
$b_{i,j}$'s is also satisfied, so Theorem 3.10 gives us a pair of isospectral
Riemannian nilmanifolds.

\medskip
\head \S4 Curvature of the Examples  \endhead 

We compare the curvature of the various examples
of isospectral manifolds constructed in sections 1-3. We continue to use the
notation
established in 1.2 and 3.1.  Since the manifolds
$G(j)$ are homogeneous, the curvature does not vary from point to
point and thus can be viewed as a tensor on the vector space ${\frak v} +
{\frak z}$
(i.e., on the Lie algebra
${\frak g}(j)$, identified with the tangent space to $G(j)$ at the identity).
The
curvatures of the manifolds $M(j)$ in 1.2 and of the closed
nilmanifolds $\njg$ in 3.1 are the same as that of $G(j)$.

The curvature of $G(j)$ is easily computed.  See \cite{E} for details.

\proclaim{Proposition 4.1}  Given inner product spaces ${\frak v}$ and ${\frak
z}$ and
a linear map $j : {\frak z} \rightarrow \so({\frak v})$, let $G(j)$ be the
associated
Riemannian manifold constructed as in 1.2.  Let $\{Z_1,\ldots, Z_r\}$ be an
orthonormal
basis of ${\frak z}$ and let $S=\frac{1}{2} \sum^r_{k=1} j^2
(Z_k)$.  For $X, Y\in {\frak v}$ and $Z, W \ in \ {\frak z}$ orthogonal
unit vectors, the sectional curvature $K$ and Ricci curvature are given as 
follows:
\item{(i)}
$$ K(X,Y) = - {\frac{3}{4}} \|[X,Y] \|^{2}
$$
$$ K(X,Z) = \frac{1}{4} \|j(Z) X \|^{2}
$$
$$ K(Z,W) = 0 
$$
\item{(ii)}
$$
\Ric (X,Y) = \langle S(X),Y \rangle 
$$
$$
\Ric (X,Z)=0
$$ 
$$ \Ric (Z,W) = -\frac{1}{4} \trace (j(Z)j(W)).
$$

In particular, if $j$ is injective, then the Ricci tensor is positive-definite
on
$\fz$ and negative semi-definite on $\fv$.

\endproclaim

\proclaim{Corollary 4.2}
Fix inner product spaces ${\frak v}$ and ${\frak z}$ and let 
$j, j^{\prime} : {\frak z} \ra \so({\frak v})$ be injective linear maps. Let
$Ric$ and $Ric^{\prime}$ denote the Ricci tensors of the associated
manifolds $G(j)$ and $G(j^{\prime})$. If $j \sim j^{\prime}$, then
$$
\Ric_{| {\frak z} \times {\frak z}} = \Ric^{\prime}_{| {\frak z}\times{\frak
z}}.
$$
\endproclaim

Thus to compare the Ricci curvatures  of the examples we need only look at 
$\Ric_{| {\frak v} \times {\frak v}}$.  The eigenvalues of $\Ric_{| {\frak v}
\times{\frak v}}$ are the eigenvalues of the operator $S$ in  Propositon 4.1.

\subhead Example 4.3\endsubhead We first consider the 7-dimensional manifolds
constructed in  Example 1.10 (see also Example 3.5). We assume that $T$ and
$T^{\prime}$ are diagonal with respect to the standard basis of ${\frak z}$
with
diagonal entries
$(a,b,c)$ and $(-a,b,c)$, respectively. Then both $\Ric_{| {\frak v}
\times {\frak v}}$ and 
$\Ric^{\prime}_{| {\frak v} \times {\frak v}}$ are diagonalized by the standard
basis of ${\frak v}$. The
four eigenvalues of $\Ric_{| {\frak v} \times {\frak v}}$ are all the
expressions of the form
${-\frac{1}{2}}\{(1 \pm a)^{2} + (1 \pm b)^{2} + (1 \pm c)^{2} \}$ with an 
even number of choices of minus signs in the terms in parentheses. The
eigenvalues
of $\Ric^{\prime}$ are obtained by changing the sign of $a$; equivalently, they
are 
all the expressions of the form above having an odd number of choices of minus 
signs. 

Thus Examples 1.10 and 3.5 yield isospectral manifolds with different Ricci
curvatures.
We note, however, that the Ricci tensors have the same norm.  
\medskip

\subhead Example 4.4\endsubhead  We consider the continuous families of
isospectral
manifolds 
$G(j_{u})$ constructed in Example 2.3, with $j_u =j_{a,b(u)}$.  Let $S_u$ be
the
operator associated with
$j_u$ as in Proposition 4.1.  We have $S_{u} = {\frac{1}{2}} (a^{2} +
b(u)^{2})$. As
noted above, the 
manifolds $G(j_{u})$, $u \in I$, have the same Ricci curvature if and only if
the 
linear operators $S_{u}, u \in I$, are isospectral.  An explicit computation
shows
that, for example, when $a$ and $b$ are chosen as in Example 3.11, then
$\det(S_{u})$ 
depends non-trivially on $u$. Thus the eigenvalues of the Ricci tensor of
$G(j_u)$
(and of $M(j_u)$) depend non-trivially on $u$. In particular, the closed
nilmanifolds
in Example 3.11 have different Ricci curvature.

However, for all choices of $a$ and $b$, $\trace(^{t}S_u S_u)$ is independent
of $u$.
Consequently, the norm of the Ricci tensor does not change during any of the 
deformations.

\medskip

 \head Appendix \endhead
\medskip

The proof of Proposition 3.7 is by an explicit calculation of the
spectra.  Using the Kirillov theory of representations of a nilpotent
Lie group, Pesce \cite{P1} computed the eigenvalues of an arbitrary
compact two-step nilmanifold.  We first summarize his results.

Let $N=(\Gamma\backslash G, g)$ be a compact two-step nilmanifold.  Thus $G$ is
a simply-connected two-step nilpotent Lie group and $g$ is a left-invariant
metric on $G$.  (We are temporarily dispensing with the notation established in
the earlier sections.)  Recall that the Laplacian of
$N$ is given by
$\Delta= -\sum_iX_i^2$, where
$\{ X_1, X_2,\ldots,X_n\}$ is an orthonormal basis of the Lie algebra $\fgs$
relative
to the inner product $<\,,\,>$ defined by $g$.  Letting $\rho =\rho_\Gamma$
denote the right action of $G$ on $L^2(N)$, then the Laplacian acts
on $L^2(N)$ as $\Delta=-\sum_i\rho_*X_i^2$.

Given any unitary representation $(V,\pi)$ of $G$ (here $V$ is a Hilbert space
and
$\pi$ is a representation of $G$ on $V$), we may define a Laplace operator
$\Delta_{g,\pi}$ on
$V$ by
$\Delta_{g,\pi}=-\sum_i\pi_*X_i^2$.  The eigenvalues of this operator
depend only on $g$ and the equivalence class of the representation
$\pi$.  The space $(L^2(N),\rho)$ is the countable direct sum of
irreducible representations $(V_\alpha,\pi_\alpha)$, each occurring
with finite multiplicity.  The spectrum of $N$ is the union, with
multiplicities,
of the spectra of the operators $\Delta_{g,\pi_{\alpha}}$.

Kirillov \cite{Ki} showed that the equivalence classes of irreducible
unitary representations of the simply-connected nilpotent Lie group
$G$ are in one to one correspondence with the orbits of the
co-adjoint action of $G$ on the dual space $\fgs^*$ of the Lie algebra
$\fgs$ of $G$.  We will denote the representation corresponding to the
co-adjoint orbit of $\sigma\in\fgs^*$ by $\pi_\sigma$.

Richardson \cite{Rn} computed the decomposition of
$L^2(\Gamma\backslash G)$
into irreducible representations $\pi_\sigma$ for an arbitrary compact
nilmanifold.  In case $G$ is two-step nilpotent, this decomposition can be
given very
explicitly.

\definition{Notation A.1}  
Given $\sigma\in \fg^*$, define
$B_{\sigma}:\fgs\times\fgs\ra R$ by
$$
B_\sigma(X,Y)=\sigma([X,Y]).
$$
Let $\fgs^\sigma =\ker (B_\sigma)$ and let $\overline{B_\sigma}$ be the
non-degenerate skew-symmetric bilinear form induced by $B_\sigma$ on
$\fgs/\fgs^\sigma$.
The image of $\log(\Gamma)$ in $\fgs/\fgs^\sigma$ is a lattice, which we
denote by ${\Cal A}_\sigma$.

We will write $\Delta_{g,\sigma}$ for $\Delta_{g,\pi_{\sigma}}$.
\enddefinition

\proclaim{Proposition A.2 (see \cite{P1})}  Let $N=(\gg,g)$ be a compact
two-step
nilmanifold, let $\fgs$ be the Lie algebra of $G$, and let $\sigma \in \fgs^*$. 
Then
$\pi_\sigma$ appears in the quasi-regular representation $\rho_\Gamma$
of $G$ on $L^2(N)$ if and only if $\sigma(\log(\Gamma)\cap\fgs^\sigma)
\subset {\bold Z}$.  In this case the multiplicity of $\pi_\sigma$ is
$m_\sigma=1$
if $\sigma(\fz)=\{0\}$, and $m_\sigma=(\det \overline{B_\sigma})^{1/2}$
otherwise, where the determinant is computed with respect to a
lattice basis of ${\Cal A}_\sigma$.
\endproclaim

Let $\bcb$ be the inner product on $\fgs^*$ defined by the Riemannian
inner product on $\fgs$.
 
\proclaim{Proposition A.3 \cite{P1}}  Let $\fz=[\fgs,\fgs]$.

\item{(a)}  If $\sigma(\fz)=0$, then $\pi_\sigma$ is a character of $G$ and 
$$
\spec(\Delta_{g,\sigma})=\{4\pi^2\|\sigma\|^2\}.
$$
\item{(b)}  If $\sigma(\fz)\neq \{0\}$, let
$\pm(-1)^{1/2}d_1,\ldots,\pm(-1)^{1/2}d_r$ be the eigenvalues of
$\overline{B_\sigma}$.  Then
$$
\spec(\Delta_{g,\sigma})=\{\mu(\sigma,p,g):p\in {\bold N}^r\}
$$
where
$$
\mu(\sigma,p,g)=4\pi^2\sum_{i=1,\ldots l} \sigma(Z_i)^2 +
2\pi\sum_{k=1,\ldots,r}(2p_k+1)d_k
$$
with $\{Z_1,\ldots, Z_l\}$ a $g$-orthonormal basis of $\fgs^\sigma$. 
The multiplicity of an eigenvalue $\mu$ is the number of $p\in {\bold N}^r$
such that $\mu =\mu(\sigma,p,g)$.

\endproclaim

\demo{Proof of Proposition 3.7} 
We use the notation of 3.1, 3.6
and A.1, and let $G=G(j)$ and $N=\njg$.  By an elementary and standard
argument, the
part of
$\spec(N)$ corresponding to all the characters $\pi_\sigma$ in part (a) of
Proposition
A.3 coincides with $\spec(N_0(j,\Gamma))$.  Thus we need only consider those
representations $\pi_\sigma$ with $\sigma(\fz)\neq 0$.

For $z$ as in Proposition 3.7 and $x,y\in\fg=\fg(j)$, observe that 
$$B_\sigma(x,y)=<[x,y],z>\sigma(z)=<j(z)x,y>\sigma(z) \leqno(1)$$  
by 1.1.  Thus $\fgs^\sigma$, as defined in A.1, coincides with ${\frak
a}+{\frak z}$. 
Hence the occurrence condition for $\pi_\sigma$ in Proposition A.2 just says
that
$\sigma_{|{\frak a}+{\frak z}}\in {\Cal K}^*$.  Observe that $\mu\in\fgs^*$
lies in the
same co-adjoint orbit as $\sigma$ if and only if $\mu_{|{\frak a}+{\frak
z}}=\sigma_{|{\frak a}+{\frak z}}$; thus we may identify co-adjoint orbits with
elements $\sigma$ of
$\kc^*$.

By equation (1), the eigenvalues of $B_\sigma$, and thus of
$\overline{B_\sigma}$, are
determined by the eigenvalues of $j(z)$ and by $\sigma(z)$.  Moreover, for
$\{Z_1\ldots Z_\ell \}$ an orthonormal basis of $\fg^\sigma={\frak a}+{\frak
z}$, we
have
$\sum_{i=1}^\ell \sigma(Z_i)^2=\|\sigma\|^2$.  Thus by Proposition A.3(b), the
eigenvalues of $\Delta_{g,\sigma}$ are completely determined by the data in
(ii) and (iii) of Proposition 3.7.

It remains to show that the data (i)--(iii) determines the multiplicity
$m_\sigma$ of $\pi_\sigma$ in the representation $\rho_\Gamma$ of $G$ on
$L^2(N)$. 
First observe that the center
$z(G)$ has Lie algebra
$\fgs^\sigma ={\frak a}+{\frak z}$.  Let $\pi:G\ra G/z(G)$ be the
projection.  The group $\pi(G)$ with the Riemannian structure induced by that
of
$G$ is Euclidean and $T:=\pi(G)/\pi(\Gamma)$ is a flat torus.  Letting
$\overline{\fgs}=\fg/\fg^\sigma$ and letting ${\Cal A}_\sigma$ be as in A.1,
then the
Lie group exponential map from
$\overline{\fgs}$ to $\pi(G)$ carries ${\Cal A}_\sigma$ to $\pi(\Gamma)$  and
induces
an isometry from the torus $\overline{\fgs}/{\Cal A}_\sigma$ to $T$, where
$\overline{\fgs}$ is given the inner product induced by that on $\fg$.  (Note
that
$\overline{\fgs}$ may be identified with the subspace ${\frak u}$ of $\fg$
defined
in 3.6.)  $N$ fibers over $T$ as a Riemannian
submersion with fiber $z(G)/(z(G)\cap\Gamma)$.  The fiber is isometric to the
torus
$({\frak a}+{\frak z})/{\Cal K}$.

Now consider the multiplicity $m_\sigma$ of $\pi_\sigma$, given in Proposition
A.2.  By equation (1), the determinant of $\overline{B_\sigma}$ with respect to
an
orthonormal basis of $\fgs/\fgs^\sigma$ (relative to the induced inner product
defined above) is determined by the eigenvalues of $j(z)$ and by $\sigma(z)$. 
To find the determinant with respect to a lattice basis of ${\Cal A}_\sigma$,
the
only additional information needed is the volume of $T$ (i.e. the ``volume''
of the lattice).  Thus it remains to show that the volume of $T$ is
determined by the data (i)--(iii).

We have two ways of viewing $N$ as a principal torus bundle over a torus;
both are Riemannian submersions.  First we have the submersion discussed
alone:
$$
\matrix z(G)/(z(G)\cap\Gamma)&\longrightarrow &N\\
                                        & &|\\
                                       & &T
\endmatrix
$$
Secondly we have a submersion with circle fiber:
$$
\matrix S:=[G,G]/([G,G]\cap\Gamma)&\hookrightarrow &N\\
                                            & &|\\
                                            & &N_0
\endmatrix
$$
The fiber circle is isometric to ${\frak z}/\kc\cap {\frak z})$ with the inner
product
$<\,,\,>$.

The second fibration and the data (i)--(iii) enable us to determine
$\Vol(N)$.  Indeed (i) gives us $\Vol(N_0)$, and from (iii) we can determine
the length of the circle $S$ as follows:  From (iii) we can find
$\min\{\|\sigma\|:\sigma\in\kc^*$ and $\|\sigma\| = |\sigma
(z)|\}=\min\{\|\sigma\|:
\sigma \in\kc^*$ and $\sigma_{|{\frak a}} =0\}$; call this
$c$.  But
$c$ is precisely the length of a basis element of the lattice in ${\frak z}^*$
dual to
$\kc\cap {\frak z}=\log(\Gamma)\cap {\frak z}$.  Hence $c$ determines the
length of the
circle $S$.  We conclude that the data (i)--(iii) determine $\Vol(N)$. 

Next the second half of the data in (iii), i\.e\., $\{\|\tau\|:\tau\in\kc^*\}$
determines the spectrum of the fiber torus in the first submersion and thus
the volume of the fiber.  This together with $\Vol(N)$ determines $\Vol(T)$.
Thus the multiplicity $m_\sigma$ is determined by the data (i)--(iii).  This
completes the proof.
\enddemo


\Refs
\medskip

\widestnumber\key{GWW}
\ref\key Be \by P. B\'{e}rard \paper Vari\'{e}t\'{e}s Riemanniennes
isospectrales nonisom\'{e}triques \jour Sem. Bourbaki \vol 705 \yr 1988-89
\endref
\ref \key B \by R. Brooks \paper Constructing isospectral manifolds \jour Amer.
Math.
Monthly \vol 95 \yr 1988 \pages 823--839 \endref
\ref \key Bd \by G. Bredon \book Introduction to compact transformation
groups \publ Academic Press\publaddr New York \yr 1972 \endref
\ref \key BT \by R. Brooks and R. Tse \paper Isospectral surfaces of small
genus
 \jour Nagoya Math. J. \vol 107 \yr1987 \pages 13--24\endref
\ref \key Bu \by P. Buser \paper  Isospectral Riemann surfaces\jour Ann. Inst.
Fourier
\vol 36 \yr 1986 \pages 167--192\endref
\ref \key D \by D. DeTurck \paper Audible and inaudible geometric properties
\jour
Rend. Sem. Fac. Sci., Univ. Cagliari \vol 58 (supplement 1988) \pages
1--26\endref
\ref \key DG \by D. DeTurck and C. Gordon \paper Isospectral deformations
II: Trace formulas, metrics, and potentials \jour Comm. Pure Appl. Math. \vol
 42 \yr 1989 \pages 1067--1095\endref
\ref\key E\by P. Eberlein
\paper Geometry of two-step nilpotent groups with a left invariant metric
\jour Ann. Scient. de l'Ecole Norm. Sup.\vol27\yr1994\pages611-660 
\endref 
\ref \key G1 \manyby C. Gordon \paper You can't hear the shape of a manifold
\inbook New Developments in Lie Theory and Their Applications \eds J. Tirao
and N. Wallace \publ Birkhauser \publaddr Boston \yr 1992\endref
\ref
\key G2
\bysame
\paper Isospectral closed Riemannian manifolds which are not locally isometric
\jour J. Diff. Geom.
\vol 37
\pages 639--649
\yr1993
\endref

\ref\key G3
\bysame
\paper Isospectral closed Riemannian manifolds which are not locally isometric,
Part II 
\inbook Contemporary Mathematics: Geometry of the Spectrum
\publ AMS\vol 173\eds R. Brooks, C. Gordon, P. Perry
\yr1994\pages 121--131
\endref
\ref\key GWW\by C. Gordon, D. Webb, and S. Wolpert
\paper Isospectral plane domains and surfaces via Riemannian orbifolds
\jour Invent. Math.
\vol 110
\pages 1--22
\yr1992
\endref
\ref\key GW1\manyby C. Gordon and E. N. Wilson
\pages 241--256
\paper Isospectral deformations of compact solvmanifolds
\yr1984
\vol 19
\jour J. Differential Geometry
\endref

\ref\key GW2\bysame
\pages 253--271
\paper The spectrum of the Laplacian on Riemannian Heisenberg manifolds
\yr1986\vol 33
\jour Mich. Math. J.
\endref

\ref\key Gt1 \manyby R. Gornet
\paper A new construction of isospectral Riemannian manifolds with examples
\jour Mich. Math. J. \vol 43 \yr 1996 \pages 159-188
\endref

\ref\key Gt2\bysame
\paper Continuous families of
Riemannian manifolds isospectral on functions but not on 1-forms\paperinfo
preprint
\endref

\ref \key I \by A. Ikeda \paper On lens spaces which are isospectral but
not isometric \jour Ann. Sci. Ec. Norm. Sup. \vol 13\yr1980 \pages
303--315\endref
\ref \key Ki \by A.A. Kirillov \paper Unitary representations of nilpotent
Lie groups\jour Russian Math. Surveys \vol 17\yr 1962\pages 53--104\endref
\ref
\key M
\by J. Milnor
\page 542
\paper Eigenvalues of the Laplace operator on certain manifolds 
\yr1964
\vol 51
\jour Proc. Nat. Acad. Sci. USA
\endref
\ref
\key P1
\manyby H. Pesce
\paper Calcul du spectre d'une nilvari\'et\'e de rang deux et applications
\vol 339
\yr1993
\pages 433--461
\jour Trans. Amer. Math. Soc.
\endref

\ref\key P2\bysame
\paper Une formule de Poisson pour les vari\'et\'es de Heisenberg
\jour Duke Math. J.\vol73\yr1994\pages79--95
\endref

\ref
\key R
\by M.S. Raghunathan
\book Discrete Subgroups of Lie Groups
\publ Springer-Verlag
\publaddr Berlin and New York
\yr1972
\endref
\ref
\key S
\by T. Sunada
\paper Riemannian coverings and isospectral manifolds
\jour Ann. of Math.
\vol 121
\yr1985
\pages 169--186
\endref
\ref
\key Sz
\by Z. Szabo
\paper Locally nonisometric yet super isospectral spaces
\paperinfo preprint
\endref
\ref
\key V
\by M. F. Vign\'eras
\pages 21--32
\paper Vari\'et\'es Riemanniennes isospectrales et non isom\'etriques
\yr1980
\vol 112
\jour Ann. of Math.
\endref
\ref
\key W
\by E. N. Wilson
\paper Isometry groups on homogeneous nilmanifolds
\jour Geom. Dedicata
\vol 12
\yr1982
\pages 337-346
\endref

\endRefs
\enddocument